\documentclass[a4, 10pt]{amsart}
\usepackage{amssymb}
\usepackage{amstext}
\usepackage{amsmath}
\usepackage{amscd}
\usepackage{latexsym}
\usepackage{amsfonts}

\theoremstyle{plain}
\newtheorem{thm}{Theorem}[section]
\newtheorem*{thma}{Theorem A}
\newtheorem*{thmb}{Theorem B}
\newtheorem*{thmf}{Theorem \ref{formula}}
\newtheorem*{thmg}{Theorem \ref{gbass}}

\newtheorem{prop}[thm]{Proposition}
\newtheorem{lem}[thm]{Lemma}

\newtheorem{cor}[thm]{Corollary}

\theoremstyle{definition}
\newtheorem{defn}[thm]{Definition}

\newtheorem{rem}[thm]{Remark}
\newtheorem{conj}[thm]{Conjecture}

\newtheorem*{quesa}{Question A}
\newtheorem*{quesb}{Question B}

\theoremstyle{remark}
\newtheorem*{pf}{{\sl Proof}}

\numberwithin{equation}{section}
\def\Hom{\mathrm{Hom}}
\def\RHom{\mathrm{{\bf R}Hom}}

\def\Ext{\mathrm{Ext}}

\def\Ker{\mathrm{Ker}}

\def\m{\mathfrak m}

\def\p{\mathfrak p}

\def\H{\mathrm{H}}
\def\Gdim{\mathrm{G}\mathrm{dim}}

\def\Gid{\mathrm{Gid}}
\def\Kdim{\mathrm{dim}}

\def\depth{\mathrm{depth}}
\def\Supp{\mathrm{Supp}}

\def\id{\mathrm{id}}

\def\grade{\mathrm{grade}}

\def\Spec{\mathrm{Spec}}

\def\G{{\mathcal G}}

\def\xx{\text{\boldmath $x$}}

\begin{document}

\title[Finite modules of finite G-injective dimension]{The existence of finitely generated modules of finite Gorenstein injective dimension}
\author{Ryo Takahashi}
\address{Department of Mathematics, Faculty of Science, Okayama University, 1-1, Naka 3-chome, Tsushima, Okayama 700-8530, Japan}
\email{takahasi@math.okayama-u.ac.jp}
\thanks{
\endgraf
}
\maketitle
\begin{abstract}
In this note, we study commutative Noetherian local rings having finitely generated modules of finite Gorenstein injective dimension.
In particular, we consider whether such rings are Cohen-Macaulay.
\end{abstract}
\section{Introduction}

Throughout this note, we assume that all rings are commutative and Noetherian.
The following celebrated theorems, where the second one had been called Bass' conjecture, are obtained by virtue of the (Peskine-Szpiro) intersection theorem.
For the details, see \cite[Corollaries 9.4.6 and 9.6.2]{BH}.

\begin{thma}
A local ring is Cohen-Macaulay if it admits a nonzero Cohen-Macaulay module of finite projective dimension.
\end{thma}

\begin{thmb}
A local ring is Cohen-Macaulay if it admits a nonzero finitely generated module of finite injective dimension.
\end{thmb}

In the sixties, Auslander \cite{Auslander} introduced {\it Gorenstein dimension} (abbr. {\it G-dimension}) as a homological invariant for finitely generated modules, and developed its notion with Bridger \cite{AB}.
This invariant is an analogue of projective dimension; any finitely generated module of finite projective dimension has finite G-dimension.

Several decades later, {\it Gorenstein projective dimension} (abbr. {\it G-projective dimension}) was defined as an extension of G-dimension to modules that are not necessarily finitely generated, and {\it Gorenstein injective dimension} (abbr. {\it G-injective dimension}) was defined as a dual version of G-projective dimension.
The notions of these dimensions are based on the work of Enochs and Jenda \cite{EJ}.
G-injective dimension is a refinement of usual injective dimension in the sense that any module of finite injective dimension has finite G-injective dimension.
So far these dimensions have been extended to complexes of modules.

Now, it is natural to ask the following questions, which are generalizations of Theorem A and Theorem B respectively.

\begin{quesa}
Is a local ring Cohen-Macaulay if it admits a nonzero Cohen-Macaulay module of finite G-dimension?
\end{quesa}

\begin{quesb}
Is a local ring Cohen-Macaulay if it admits a nonzero finitely generated module of finite G-injective dimension?
\end{quesb}

Both of these questions are presented in Christensen's book \cite{Christensen}, while the second one is also found in \cite{KY}.

The author proved that Question A has an affirmative answer for any local ring of type one:

\begin{thm}\cite[Theorem 2.3]{Takahashi}\label{hungar}
The following are equivalent for a local ring $R$:
\begin{enumerate}
\item[$(1)$]
$R$ is Gorenstein;
\item[$(2)$]
$R$ admits an ideal $I$ of finite G-dimension such that $R/I$ is Gorenstein;
\item[$(3)$]
$R$ admits a nonzero Cohen-Macaulay module of type one and of finite G-dimension;
\item[$(4)$]
$R$ has type one and admits a nonzero Cohen-Macaulay module of finite G-dimension.
\end{enumerate}
\end{thm}

In this note, we will consider Question B.
The main theorems of this note are the following two theorems.

\begin{thmf}
Let $R$ be a local ring with a dualizing complex.
Then for a nonzero finitely generated $R$-module $M$, there is an inequality
$$
\Kdim\, M \leq \Gid\, M.
$$
\end{thmf}

\begin{thmg}
Let $R$ be a local ring with a dualizing complex, and let $M$ be a nonzero finitely generated $R$-module of finite G-injective dimension.
\begin{enumerate}
\item[$(1)$]
If $\Kdim\, M = \Kdim\, R$, then $R$ is Cohen-Macaulay.
\item[$(2)$]
\begin{enumerate}
\item[{\rm (i)}]
If $M$ is cyclic, then $R$ has type one.
\item[{\rm (ii)}]
Suppose that $R$ has type one.
If $M$ is Cohen-Macaulay or $\depth\, M\geq\depth\, R$, then $R$ is Gorenstein.
\end{enumerate}
\end{enumerate}
\end{thmg}

In the next section we will prepare some definitions and lemmas, and in the third (final) section we will state and prove our results.

\section{Preliminaries}

In this section, we will write down only results to use later.
Throughout this section, let $R$ be a local ring.
In the following, an {\it $R$-complex} will mean a chain complex $X=(\cdots \to X_{i+1} \to X_i \to X_{i-1} \to \cdots )$ of $R$-modules, and a {\it bounded} $R$-complex will mean a homologically bounded $R$-complex, namely an $R$-complex $X$ such that $\H _i(X)=0$ for any $i\gg 0$ and $i\ll 0$.
An {\it isomorphism} of $R$-complexes will mean an isomorphism in the derived category of the category of $R$-modules.

\begin{defn}\cite[(2.3.2)]{Christensen}
\begin{enumerate}
\item[$(1)$] We denote by $\G (R)$ the class of all finitely generated $R$-modules $M$ satisfying the following three conditions:
\begin{enumerate}
\item[{\rm (i)}] The natural homomorphism $M\to\Hom _R (\Hom _R (M, R), R)$ is an isomorphism;
\item[{\rm (ii)}] $\Ext _R ^i (M, R)=0$ for any $i>0$;
\item[{\rm (iii)}] $\Ext _R ^i (\Hom _R (M, R), R)=0$ for any $i>0$.
\end{enumerate}
\item[$(2)$] Let $X$ be an $R$-complex and $n$ an integer.
If there exists an $R$-complex of this form
$$
A=(0 \to A_n \to A_{n-1} \to \cdots \to A_{m+1} \to A_m \to 0)
$$
which is isomorphic to $X$ such that $A_i\in\G (R)$ for any $i$, then we say that $X$ has {\it G-dimension at most $n$} and write $\Gdim _R\, X\leq n$.
\end{enumerate}
\end{defn}

Let $X$ be an $R$-complex.
If $X\cong 0$ i.e. $X$ is exact, then we set $\Gdim _R\, X=-\infty$.
If $\Gdim _R\, X\leq n$ for some integer $n$, then we say that $X$ has finite G-dimension and write $\Gdim _R\, X<\infty$.
If $\Gdim _R\, X\not\leq n$ for any integer $n$, then we say that $X$ has infinite G-dimension and write $\Gdim _R\, X=\infty$.
If $\Gdim _R\, X\leq n$ but $\Gdim _R\, X\not\leq n-1$, then we say that $X$ has G-dimension $n$ and write $\Gdim _R\, X=n$.

\begin{lem}\cite[(2.3.8)(2.3.13)]{Christensen}\cite[8.7]{Avramov}\label{ref}
\begin{enumerate}
\item[$(1)$] The following are equivalent for a bounded $R$-complex $X$ with finitely generated homology:
\begin{enumerate}
\item[{\rm (i)}] $\Gdim _R\, X<\infty$;
\item[{\rm (ii)}] The complex $\RHom _R (X, R)$ is bounded, and the natural homomorphism $X \to \RHom _R (\RHom _R (X, R), R)$ is an isomorphism.
\end{enumerate}
\item[$(2)$] Let $X$ be a bounded $R$-complex with finitely generated homology.
If $\Gdim _R\, X<\infty$, then
$$
\Gdim _R\, X = -\inf\RHom _R (X, R) = \depth\, R - \depth _R\, X.
$$
\item[$(3)$] Let $M$ be a finitely generated $R$-module.
\begin{enumerate}
\item[{\rm (i)}] For an $M$-regular element $x\in R$,
$$
\Gdim _R\, M/xM=\Gdim _R\, M+1.
$$
\item[{\rm (ii)}] Let $R\to S$ be a faithfully flat homomorphism of local rings.
Then
$$
\Gdim _S\, M\otimes _R S = \Gdim _R\, M.
$$
\end{enumerate}
\end{enumerate}
\end{lem}

\begin{defn}\cite[(6.1.1)(6.2.2)]{Christensen}
\begin{enumerate}
\item[$(1)$] An $R$-module $M$ is said to be {\it G-injective} if there exists an exact complex
$$
C=(\cdots \overset{d_3}{\longrightarrow} C_2 \overset{d_2}{\longrightarrow} C_1 \overset{d_1}{\longrightarrow} C_0 \overset{d_0}{\longrightarrow} C_{-1} \overset{d_{-1}}{\longrightarrow} \cdots)
$$
of injective $R$-modules with $M\cong \Ker (d_0)$ such that the complex $\Hom _R (T, C)$ is also exact for any injective $R$-module $T$.
\item[$(2)$] Let $X$ be an $R$-complex and $n$ an integer.
If there exists an $R$-complex of this form
$$
A=(0 \to A_m \to A_{m-1} \to \cdots \to A_{-n+1} \to A_{-n} \to 0)
$$
which is isomorphic to $X$ such that $A_i$ is G-injective for any $i$, then we say that $X$ has {\it G-injective dimension at most $n$} and write $\Gid _R\, X\leq n$.
\end{enumerate}
\end{defn}

The conditions $\Gid _R\, X=-\infty$, $\Gid _R\, X<\infty$, $\Gid _R\, X=\infty$, and $\Gid _R\, X=n$ are defined similarly to G-dimension.

\begin{lem}\cite[(6.2.3)]{Christensen}\cite[(2.12)(6.3)(6.4)]{CFH}\label{gid}
\begin{enumerate}
\item[$(1)$]
For a bounded $R$-complex $X$,
$$
\Gid _R\, X \leq \id _R\, X.
$$
\item[$(2)$]
Assume that $R$ admits a dualizing complex $D$.
\begin{enumerate}
\item[{\rm (i)}]
For a bounded $R$-complex $X$ with finitely generated homology,
$$
\Gid _R\, X <\infty \ \Leftrightarrow\  \Gdim _R\, \RHom _R (X, D) <\infty.
$$
\item[{\rm (ii)}]
Let $M$ be a nonzero finitely generated $R$-module of finite G-injective dimension.
Then
$$
\Gid _R\, M=\depth\, R.
$$
\end{enumerate}
\end{enumerate}
\end{lem}

\section{Main theorems}

First of all, we remark that Theorem B implies Theorem A.
Indeed, let $R$ be a local ring and $M$ a Cohen-Macaulay $R$-module of finite projective dimension.
Setting $\depth\,M=\Kdim\,M=n$, we have an $M$-sequence $\xx = x_1, x_2, \dots, x_n$ in $R$ such that $M/\xx M$ has finite length.
Then since the $R$-module $M/\xx M$ also has finite projective dimension, the dual $R$-module $(M/\xx M)^{\vee}$ has finite length and finite injective dimension, where $(-)^{\vee}$ denotes the Matlis dual over $R$.
Hence it follows from Theorem B that $R$ is Cohen-Macaulay.
Thus Theorem B implies Theorem A.

Similarly, we have the following.

\begin{prop}
If Question B has an affirmative answer, then so has Question A.
\end{prop}

\begin{pf}
Let $R$ be a local ring and $M$ a Cohen-Macaulay $R$-module of finite G-dimension.
We want to show that $R$ is Cohen-Macaulay under the assumption that Question B has an affirmative answer.
Replacing $R$ and $M$ by their completions respectively, we may assume that $R$ is complete; see Lemma \ref{ref}(3)(ii).
Since $M$ is Cohen-Macaulay $R$-module, there exists an $M$-sequence $\xx$ in $R$ such that $M/\xx M$ has finite length.
Lemma \ref{ref}(3)(i) says that the $R$-module $M/\xx M$ has finite G-dimension.
Replacing $M$ by $M/\xx M$, we may assume that $M$ has finite length and finite G-dimension.
Note that $\RHom _R (M, D)$ is isomorphic to $M^{\vee}$ up to shift, where $D$ is a dualizing complex of $R$.
Hence, according to Lemma \ref{gid}(2)(i), $M^{\vee}$ has finite length and finite G-injective dimension.
Since we are assuming that Question B has an affirmative answer, we conclude that $R$ is Cohen-Macaulay, as desired.
\qed
\end{pf}

From now on, we will study Question B over a local ring admitting a dualizing complex, and prove the main theorems of this note.

Let $(R, \m, k)$ be a local ring.
We set $d=\Kdim\, R$ and $t=\depth\, R$.
Let $D=(0 \to D_0 \to D_{-1} \to \cdots \to D_{-d} \to 0)$ be a dualizing complex of $R$, where $D_i = \bigoplus _{\p\in\Spec\, R,\ \Kdim\, R/\p =d+i} E_R (R/\p )$.
The local duality theorem yields $\inf D = t-d$.

Let $M\neq 0$ be a finitely generated $R$-module with $\Gid\,M<\infty$.
Put $X=\RHom _R (M, D)$.
Lemma \ref{gid}(2)(i) implies $\Gdim\,X<\infty$.
Since $\depth\,X = \depth\,D = d$ by \cite[(A.6.4)]{Christensen}, we have $\inf\RHom _R (X, R) = d-t$ by Lemma \ref{ref}(2).
Since $\RHom _R (X, R) \cong \RHom _R (D, M)$, we have $\sup\RHom _R (X, R) = \sup\RHom _R (D, M) \leq d-t$ by \cite[(A.4.6.1)]{Christensen}.
Hence $\sup\RHom _R (X, R) = \inf\RHom _R (X, R) = d-t$.
Setting $N=\H _{d-t}(\RHom _R (X, R))$, we get $\RHom _R (X, R) \cong N[d-t]$.
Lemma \ref{ref}(1) implies
\begin{equation}\label{gdimn}
\Gdim\, N < \infty.
\end{equation}
Noting that $N \cong \H _{d-t} (\RHom _R (D, M))$ and $\inf D = t-d$, we obtain
\begin{equation}\label{shape}
N \cong \Hom _R (\H _{t-d}(D), M).
\end{equation}

\begin{rem}
In the above part, we constructed the complex $\RHom _R (D, M)$ of finite G-dimension from a module $M$ of finite G-injective dimension.
A deeper investigation in connection with this fact has been made by Foxby; he has given an equivalence between the Auslander class and the Bass class, which are certain full subcategories of the category of complexes.
This equivalence is said to be {\it Foxby equivalence}.
For the definition and the properties, see \cite[Section 3.3]{Christensen}.
\end{rem}

\begin{lem}\label{ineq}
With the notation introduced above,
$$
\depth\, M = \depth\, N \leq \Kdim\, N \leq \Kdim\, M \leq t.
$$
\end{lem}

\begin{pf}
By \cite[(A.8.5.1), (A.8.5.2)]{Christensen} we get $\inf X=\depth\, M-d$ and $\sup X=\Kdim\, M-d$.
Note that $\RHom _R (N, R) \cong X[d-t]$.
By \eqref{gdimn} and Lemma \ref{ref}(2) we obtain $t-\depth\, N = \Gdim\, N = -\inf\RHom _R (N, R) = -\inf X - (d-t) = t- \depth\, M$, and $0 \leq \grade\, N = -\sup\RHom _R (N, R) = -\sup X - (d-t) = t-\Kdim\, M$.
Hence $\depth\, N =\depth\, M$ and $\Kdim\, M \leq t$.
On the other hand, since $\Supp\, N \subseteq \Supp\, M$ by \eqref{shape}, we have $\Kdim\, N \leq \Kdim\, M$.
Thus the proof of the lemma is completed.
\qed
\end{pf}

Let $R$ be a local ring and $M$ a nonzero finitely generated $R$-module.
Then it is well-known that there is an inequality
$$
\Kdim\, M \leq \id\, M.
$$
(For the proof, see \cite[Theorem 3.1.17]{BH}.)
Combining Lemma \ref{ineq} with Lemma \ref{gid}(2)(ii) extends this formula of injective dimension to that of G-injective dimension.

\begin{thm}\label{formula}
Let $R$ be a local ring with a dualizing complex.
Then for a nonzero finitely generated $R$-module $M$, there is an inequality
$$
\Kdim\, M \leq \Gid\, M.
$$
\end{thm}

Recall that the {\it type} of a finitely generated module $M$ over a local ring $R$ with residue field $k$ is defined to be the dimension of the $k$-vector space $\Ext _R ^t (k, M)$ where $t=\depth _R\, M$.
Using Lemma \ref{ineq} again, we can give a result with relation to Question B in the first section of this note:

\begin{thm}\label{gbass}
Let $R$ be a local ring with a dualizing complex, and let $M$ be a nonzero finitely generated $R$-module of finite G-injective dimension.
\begin{enumerate}
\item[$(1)$]
If $\Kdim\, M = \Kdim\, R$, then $R$ is Cohen-Macaulay.
\item[$(2)$]
\begin{enumerate}
\item[{\rm (i)}]
If $M$ is cyclic, then $R$ has type one.
\item[{\rm (ii)}]
Suppose that $R$ has type one.
If $M$ is Cohen-Macaulay or $\depth\, M\geq\depth\, R$, then $R$ is Gorenstein.
\end{enumerate}
\end{enumerate}
\end{thm}

\begin{pf}
(1) This assertion immediately follows from Lemma \ref{ineq}.

(2) We use the same notation as in the first part of this section.

(i) Set $Y=\RHom _R (X, R)$.
Using Lemma \ref{ref}(1), we get isomorphisms
$$
\RHom _R (M, k)[-d] \cong \RHom _R (k, X) \cong \RHom _R (Y, \RHom _R (k, R)).
$$
Since $M$ is cyclic, comparing the $(-d)$th homology modules, we obtain $k \cong \Ext _R ^t (k, R)$ (cf. \cite[(A.4.6)]{Christensen}).
This says that $R$ has type one.

(ii) It follows from Lemma \ref{ineq} and \eqref{gdimn} that $N$ is a Cohen-Macaulay $R$-module of finite G-dimension.
Since $R$ has type one, Theorem \ref{hungar} implies that $R$ is Gorenstein.
\qed
\end{pf}

\begin{cor}
Let $R$ be a local ring of type one with a dualizing complex.
If there is a finitely generated G-injective $R$-module, then $R$ is Artinian Gorenstein.
\end{cor}

\begin{pf}
Let $M$ be a finitely generated G-injective $R$-module.
Since $\Gid\, M=0$, we have $\depth\, R=0$ by Lemma \ref{gid}(2)(ii).
Therefore Theorem \ref{gbass}(2)(ii) implies that $R$ is Gorenstein, and $\Kdim\, R=\depth\, R=0$.
\qed
\end{pf}

Holm \cite[Theorem 2.1]{Holm} proved that any ring of finite G-injective dimension is Gorenstein.
Letting $M=R$ in Theorem \ref{gbass}(2), we get his theorem for a local ring with a dualizing complex:

\begin{cor}[Holm]
Let $R$ be a local ring with a dualizing complex.
If $\Gid\, R<\infty$, then $R$ is Gorenstein.
\end{cor}

Peskine and Szpiro \cite[Chapitre II, Th\'{e}or\`{e}me (5.5)]{PS} proved that a local ring is Gorenstein if it has a nonzero cyclic module of finite injective dimension.
Taking Lemma \ref{gid}(1) and the above results into account, we end this note by stating a conjecture.

\begin{conj}
Let $R$ be a local ring.
If there exists a nonzero cyclic $R$-module $M$ of finite G-injective dimension, then $R$ is Gorenstein.
\end{conj}

{\sc Acknowledgments.}
The author would like to thank an anonymous referee for careful reading and helpful comments.


\end{document}